\documentclass[12pt]{article} 
\usepackage{amsfonts,latexsym}
\usepackage{amsthm}
\input{psfig.sty}

\newtheorem{theorem}{Theorem}[section]
\newtheorem{proposition}[theorem]{Proposition}
\newtheorem{lemma}[theorem]{Lemma}

\theoremstyle{definition}

\theoremstyle{remark}

\begin{document}
\newcommand{\beq}{\begin{equation}} \newcommand{\eeq}{\end{equation}}
\newcommand{\beas}{\begin{eqnarray*}}
\newcommand{\eeas}{\end{eqnarray*}} 
\newcommand{\zz}{\mathbb{Z}}
\newcommand{\pp}{\mathbb{P}} 
\newcommand{\nn}{\mathbb{N}}
\newcommand{\rr}{\mathbb{R}}
\newcommand{\bm}[1]{{\mbox{\boldmath $#1$}}}
\newcommand{\con}{\mathrm{Comp}(n)}
\newcommand{\sn}{\mathfrak{S}_n} 
\newcommand{\fs}{\mathfrak{S}}
\newcommand{\st}{\,:\,} 

\thispagestyle{empty}

\vskip 20pt
\begin{center}
{\large\bf The Descent Set and Connectivity Set of a
  Permutation}
\vskip 15pt
{\bf Richard P. Stanley}\footnote{Partially supported by
  NSF grant \#DMS-9988459 and by the Institut Mittag-Leffler.}\\
{\it Department of Mathematics, Massachusetts Institute of
Technology}\\
{\it Cambridge, MA 02139, USA}\\
{\texttt{rstan@math.mit.edu}}\\[.2in]
{\bf\small version of 11 July 2005}\\
\end{center}

\begin{abstract}
The descent set $D(w)$ of a permutation $w$ of $1,2,\dots,n$ is a
standard and well-studied statistic. We introduce a new statistic, the
\emph{connectivity set} $C(w)$, and show that it is a kind of dual
object to $D(w)$. The duality is stated in terms of the inverse of a
matrix that records the joint distribution of $D(w)$ and $C(w)$. We
also give a variation involving permutations of a multiset and a
$q$-analogue that keeps track of the number of inversions of $w$.
\end{abstract}

\section{A duality between descents and connectivity.} \label{sec1}
\indent Let $\sn$ denote the symmetric group of permutations of
$[n]=\{1,2,\dots,n\}$, and let $w=a_1a_2\cdots a_n\in\sn$. The
\emph{descent set} $D(w)$ is defined by
  $$ D(w) = \{i\st a_i>a_{i+1}\}\subseteq [n-1]. $$
The descent set is a well-known and much studied statistic on
permutations with many applications, e.g.,
\cite[Exam.~2.24, Thm.~3.12.1]{ec1}\cite[{\S}7.23]{ec2}. Now define
the \emph{connectivity set} $C(w)$ by
  \beq C(w) = \{i\st a_j<a_k\ \mathrm{for\ all}\ j\leq i<k\}. 
    \label{eq:cwdef} \eeq
The connectivity set seems not to have been considered before except
for equivalent definitions by Comtet \cite[Exer.~VI.14]{comtet2} and
Callan \cite{callan}
with no further development, but some related notions have been
investigated. In particular, a permutation $w$ with $C(w)=\emptyset$
is called \emph{connected} or \emph{indecomposable}. If $f(n)$ denotes
the number of connected permutations, then Comtet
\cite[Exer.~VI.14]{comtet2} showed that
  $$ \sum_{n\geq 1} f(n)x^n = 1-\frac{1}{\sum_{n\geq 0} n!x^n}, $$ 
and he also considered the number $\#C(w)$ of components. For further
references on connected permutations, see \cite{sloane}. In this paper
we will establish a kind of ``duality'' between descent sets and
connectivity sets.

We write $S=\{i_1,\dots,i_k\}_<$ to denote that $S=\{i_1,\dots,i_k\}$
and $i_1<\cdots<i_k$. Given $S=\{i_1,\dots,i_k\}_<\subseteq [n-1]$, 
define 
  $$ \eta(S) = i_1!\,(i_2-i_1)!\cdots (i_k-i_{k-1})!\,
  (n-i_k)!. $$
Note that $\eta(S)$ depends not only on $S$ but also on $n$. The
integer $n$ will always be clear from the context.  The first
indication of a duality between $C$ and $D$ is the following result.

\begin{proposition}
Let $S\subseteq [n-1]$. Then
  \beas \#\{w\in\sn\st S\subseteq C(w)\} & = & \eta(S)\\
        \#\{w\in\sn\st S\supseteq D(w)\} & = &
         \frac{n!}{\eta(S)}. \eeas
\end{proposition}

\textbf{Proof.} The result for $D(w)$ is well-known, e.g.,
\cite[Prop.~1.3.11]{ec1}. To obtain a permutation $w$ satisfying
$S\supseteq D(w)$, choose an ordered partition $(A_1,\dots,A_{k+1})$
of $[n]$ with $\#A_j=i_j-i_{j-1}$ (with $i_0=0$, $i_{k+1}=n$) in
$n!/\eta(S)$ ways, then arrange the elements of $A_1$ in increasing
order, followed by the elements of $A_2$ in increasing order, etc.

Similarly, to obtain a permutation $w$ satisfying $S\subseteq C(w)$,
choose a permutation of $[i_1]$ in $i_1!$ ways, followed by a
permutation of $[i_1+1,i_2]:=\{i_1+1,i_1+2,\dots,i_2\}$ in
$(i_2-i_1)!$ ways, etc. $\ \Box$ 

\medskip
Let $S,T\subseteq [n-1]$. Our main interest is in the joint
distribution of the statistics $C$ and $D$, i.e., in the numbers
  $$ \Gamma_{ST} = \#\{w\in\sn\st C(w)=\overline{S},\ D(w)=T\}, $$
where $\overline{S}=[n-1]-S$. (It will be more notationally convenient
to use this definition of $\Gamma_{ST}$ rather than having $C(w)=S$.)
To this end, define
  \begin{eqnarray} A_{ST} & = & \#\{w\in\sn\st \overline{S}\subseteq
    C(w),\ T\subseteq 
    D(w)\} \nonumber\\ & = & \sum_{{S'\supseteq S\atop T'\supseteq
    T}}\Gamma_{S'T'}. \label{eq:ab} \end{eqnarray}
For instance, if $n=4$, $S=\{2,3\}$, and $T=\{3\}$, then $A_{ST}=3$,
corresponding to the permutations 1243, 1342, 1432, while
$\Gamma_{ST}=1$, corresponding to 1342. Tables of $\Gamma_{ST}$ for
$n=3$ and $n=4$ are given in Figure~\ref{fig:t1}, and for $n=5$ in
Figure~\ref{fig:t2}.

\begin{figure}
\begin{center}
\begin{tabular}{c|cccc}
  $S\backslash T$ & $\emptyset$ & 1 & 2 & 12\\ \hline
      $\emptyset$ &  1 \\
                1 &  0 & 1\\
                2 &  0 & 0 & 1\\
               12 &  0 & 1 & 1 & 1\\
\end{tabular}

\medskip
\begin{tabular}{c|cccccccc}
  $S\backslash T$ & $\emptyset$ & 1 & 2 & 3 & 12 & 13 & 23
     & 123\\ 
     \hline
 $\emptyset$ & 1\\
           1 & 0 & 1\\
           2 & 0 & 0 & 1\\
           3 & 0 & 0 & 0 & 1\\
          12 & 0 & 1 & 1 & 0 & 1\\
          13 & 0 & 0 & 0 & 0 & 0 & 1\\
          23 & 0 & 0 & 1 & 1 & 0 & 0 & 1\\
         123 & 0 & 1 & 2 & 1 & 2 & 4 & 2 & 1\\
\end{tabular}       
\end{center}   
\caption{Table of $\Gamma_{ST}$ for $n=3$ and $n=4$}
\label{fig:t1}
\end{figure}

\begin{figure}
\begin{center}
\begin{tabular}{c|cccccccccccccccc}
  $S\backslash T$ & $\emptyset$ & 1 & 2 & 3 & 4 & 12 & 13 &
   14 & 23 & 24 & 34 & 123 & 124 & 134 & 234 & 1234\\ \hline
  $\emptyset$ & 1\\
     1 & 0 & 1\\
     2 & 0 & 0 & 1\\
     3 & 0 & 0 & 0 & 1\\
     4 & 0 & 0 & 0 & 0 & 1\\
    12 & 0 & 1 & 1 & 0 & 0 & 1\\
    13 & 0 & 0 & 0 & 0 & 0 & 0 & 1\\
    14 & 0 & 0 & 0 & 0 & 0 & 0 & 0 & 1\\
    23 & 0 & 0 & 1 & 1 & 0 & 0 & 0 & 0 & 1\\
    24 & 0 & 0 & 0 & 0 & 0 & 0 & 0 & 0 & 0 & 1\\
    34 & 0 & 0 & 0 & 1 & 1 & 0 & 0 & 0 & 0 & 0 & 1\\
   123 & 0 & 1 & 2 & 1 & 0 & 2 & 4 & 0 & 2 & 0 & 0 & 1\\
   124 & 0 & 0 & 0 & 0 & 0 & 0 & 0 & 1 & 0 & 1 & 0 & 0 & 1\\
   134 & 0 & 0 & 0 & 0 & 0 & 0 & 1 & 1 & 0 & 0 & 0 & 0 & 0 & 1\\
   234 & 0 & 0 & 1 & 2 & 1 & 0 & 0 & 0 & 2 & 4 & 2 & 0 & 0 & 0 & 1\\
  1234 & 0 & 1 & 3 & 3 & 1 & 3 & 10 & 8 & 6 & 10 & 3 & 3 & 8 & 
     8 & 3 & 1\\   
\end{tabular}
\end{center}
\caption{Table of $\Gamma_{ST}$ for $n=5$}
\label{fig:t2}
\end{figure}

\begin{theorem} \label{thm:alpha}
We have
  $$ A_{ST}=\left\{ \begin{array}{rl}
        \eta(\overline{S})/\eta(\overline{T}), &   
        \mathrm{if}\ \overline{S}\cap T=\emptyset;\\ 0, &
        \mathrm{otherwise}, \end{array} \right. $$
\end{theorem}

\textbf{Proof.} Let $w=a_1\cdots a_n\in\sn$. If $i\in C(w)$ then
$a_i<a_{i+1}$, so $i\not\in D(w)$. Hence $A_{ST}=0$ if
$\overline{S}\cap T\neq\emptyset$.

Assume therefore that $\overline{S}\cap T=\emptyset$. Let
$C(w)=\{c_1,\dots,c_j\}_<$ with $c_0=0$ and $c_{j+1}=n$. Fix $0\leq
h\leq j$, and let
  $$ [c_h,c_{h+1}]\cap \overline{T}= \{
     c_h=i_1,i_2,\dots,i_k=c_{h+1}\}_<. $$
If $w=a_1\cdots a_n$ with $\overline{S}\subseteq C(w)$ and $T\subseteq
D(w)$, then the number of choices for $a_{c_h}+1,
a_{c_h}+2,\dots,a_{c_{h+1}}$ is just the multinomial coefficient 
  $$ {c_{h+1}-c_h\choose i_2-i_1, i_3-i_2, \dots,i_k-i_{k-1}}:=
    \displaystyle\frac{(c_{h+1}-c_h)!}{(i_2-i_1)!\,(i_3-i_2)!\cdots
     (i_k-i_{k-1})!}. $$ 
Taking the product over all $0\leq h\leq j$ yields
$\eta(\overline{S})/\eta(\overline{T})$. $\ \Box$ 

\medskip
Theorem~\ref{thm:alpha} can be restated matrix-theoretically. Let $M=
(M_{ST})$ be the matrix whose rows and columns are indexed by subsets
$S,T\subseteq [n-1]$ (taken in some order), with 
  $$ M_{ST}=\left\{ \begin{array}{rl} 1, & \mathrm{if}\ S\supseteq T;
    \\ 0, & \mathrm{otherwise}. \end{array} \right. $$
Let $D=(D_{ST})$ be the diagonal matrix with $D_{SS}=
\eta(\overline{S})$. Let $A=(A_{ST})$, i.e., the matrix
whose $(S,T)$-entry is $A_{ST}$ as defined in
(\ref{eq:ab}). Then it is straightforward to 
check that Theorem~\ref{thm:alpha} can be restated as follows:
  \beq A = DMD^{-1}. \label{eq:admd} \eeq
Similarly, let $\Gamma=(\Gamma_{ST})$. Then it is immediate
from equations~(\ref{eq:ab}) and (\ref{eq:admd}) that
  \beq M\Gamma M = A. \label{eq:mbm} \eeq
\indent The main result of this section (Theorem~\ref{thm:main} below)
computes the inverse of the matrices $A$, $\Gamma$, and a matrix $B$
intermediate between $A$ and $\Gamma$. The matrix $B$ arose from the
theory of quasisymmetric functions in response to a question from
Louis Billera and Vic Reiner and was the original motivation for this
paper, as explained in the Note below. See for example
\cite[{\S}7.19]{ec2} for an introduction to quasisymmetric functions.
We will not use quasisymmetric functions elsewhere in this paper.

Let $\con$ denote the set of all compositions $\alpha=
(\alpha_1,\dots,\alpha_k)$ of $n$, i.e, $\alpha_i\geq 1$ and $\sum
\alpha_i=n$.  Let $\alpha= (\alpha_1,\dots,\alpha_k) \in\con$, and let
$\fs_\alpha$ denote the subgroup of $\sn$ consisting of all
permutations $w=a_1\cdots a_n$ such that $\{1,\dots,\alpha_1\}=
\{a_1,\dots,a_{\alpha_1}\}$, $\{\alpha_1+1,\dots, \alpha_1+\alpha_2\}=
\{a_{\alpha_1+1},\dots,a_{\alpha_1 +\alpha_2}\}$, etc. Thus
$\fs_\alpha\cong \fs_{\alpha_1}\times\cdots\times\fs_{\alpha_k}$ and
$\#\fs_\alpha=\eta(S)$, where $S=\{\alpha_1,\alpha_1+\alpha_2,\dots,
\alpha_1+\cdots+\alpha_{k-1}\}$. If $w\in\sn$ and $D(w)=
\{i_1,\dots,i_k\}_<$, then define the \emph{descent composition}
co$(w)$ by
  $$ \mathrm{co}(w)=(i_1,i_2-i_1,\dots,i_k-i_{k-1},n-i_k)\in\con. $$
Let $L_\alpha$ denote the fundamental quasisymmetric function indexed
by $\alpha$ \cite[(7.89)]{ec2}, and define 
  \beq R_\alpha = \sum_{w\in \fs_\alpha} L_{\mathrm{co}(w)}. 
     \label{eq:palpha} \eeq
It is easy to see that the set $\{R_\alpha\,\st\,\alpha\in\con\}$ is a
$\zz$-basis for the additive group of all homogeneous quasisymmetric
functions over $\zz$ of degree $n$. In fact, the transition matrix
between the bases $L_\alpha$ and $R_\alpha$ is lower unitriangular
(with respect to a suitable ordering of the rows and columns), as is
immediate from equation (\ref{eq:bst}) below. 

Given $\alpha=(\alpha_1,\dots,\alpha_k)\in\con$, let
$S_\alpha=\{\alpha_1,\alpha_1+\alpha_2,\dots,\alpha_1+
\cdots+\alpha_{k-1}\}$. 
Note that $w\in\fs_\alpha$ if and only if $S_\alpha\subseteq
C(w)$. Hence equation (\ref{eq:palpha}) can be rewritten as
  $$ R_\alpha = \sum_\beta B_{\overline{S_\alpha} S_\beta}L_\beta, $$
where
  \beq B_{ST} = \#\{w\in\sn\st \overline{S}\subseteq C(w),\ T=D(w)\}. 
    \label{eq:bst} \eeq
In particular, the problem of expressing the $L_\beta$'s as linear
combinations of the $R_\alpha$'s is equivalent to inverting the matrix
$B=(B_{ST})$.

\medskip
\textsc{Note.} The question of Billera and Reiner mentioned above is
the following.  Let $P$ be a finite poset, and define the
quasisymmetric function
  $$ K_P = \sum_f x^f, $$
where $f$ ranges over all order-preserving maps $f:P\rightarrow
\{1,2,\dots\}$ and $x^f = \prod_{t\in P}x_{f(t)}$ (see
\cite[(7.92)]{ec2}). Billera and Reiner asked whether the
quasisymmetric functions $K_P$ generate (as a $\zz$-algebra) or even
span (as an additive abelian group) the space of all quasisymmetric
functions. Let $\bm{m}$ denote an $m$-element antichain. The
\emph{ordinal sum} $P\oplus Q$ of two posets $P,Q$ with disjoint
elements is the poset on the union of their elements satisfying $s\leq
t$ if either (1) $s,t\in P$, (2) $s,t\in Q$ or (3) $s\in P$ and $t\in
Q$. If $\alpha=(\alpha_1,\dots,\alpha_k)\in\con$ then let
$P_\alpha=\bm{\alpha_1}\oplus\cdots \oplus\bm{\alpha_k}$. It is easy to
see that $K_{P_\alpha}=R_\alpha$, so the $K_{P_\alpha}$'s form a
$\zz$-basis for the homogeneous quasisymmetric functions of degree $n$,
thereby answering the question of Billera and Reiner.

\medskip
It is immediate from the definition of matrix multiplication and
(\ref{eq:mbm}) that the matrix $B$ satisfies
  \beq B = M\Gamma = AM^{-1}. \label{eq:bmc} \eeq
In view of equations~(\ref{eq:admd}), (\ref{eq:mbm}) and
(\ref{eq:bmc}) the computation of $A^{-1}$, $B^{-1}$, and $\Gamma^{-1}$
will reduce to computing $M^{-1}$, which is a simple and well-known
result. For any invertible matrix $N=(N_{ST})$, write $N^{-1}_{ST}$
for the $(S,T)$-entry of $N^{-1}$.

\begin{lemma} \label{lemma:minv}
We have 
  \beq M^{-1}_{ST} = (-1)^{\#S+\#T}M_{ST}. \label{eq:minv} \eeq
\end{lemma}

\textbf{Proof.} Let $f,g$ be functions from subsets of $[n]$ to $\rr$
(say) related by
  \beq f(S) = \sum_{T\subseteq S} g(T). \label{eq:ie1} \eeq
Equation~(\ref{eq:minv}) is then equivalent to the inversion formula
  \beq g(S) = \sum_{T\subseteq S}(-1)^{\#(S-T)}f(T). \label{eq:ie2} \eeq
This is a standard combinatorial result with many proofs, e.g.,
\cite[Thm.~2.1.1, Exam.~3.8.3]{ec1}. $\ \Box$

\textsc{Note.} The matrix $M$ represents the zeta function of the
boolean algebra ${\cal B}_n$ \cite[{\S}3.6]{ec1}. Hence
Lemma~\ref{lemma:minv} can be regarded as the determination of the
M\"obius function of ${\cal B}_n$ \cite[Exam.~3.8.3]{ec1}. All our
results can easily be formulated in terms of the incidence algebra of
${\cal B}_n$.

\begin{theorem} \label{thm:main}
The matrices $A,B,\Gamma$ have the following inverses:
  \begin{eqnarray} A^{-1}_{ST} & = & (-1)^{\#S+\#T}A_{ST}
    \label{eq:aist} \\ 
        B^{-1}_{ST} & = & (-1)^{\#S+\#T}\#\{ w\in\sn\st
          \overline{S}= C(w),\ T\subseteq D(w)\} \label{eq:bist}\\
        \Gamma^{-1}_{ST} & = & (-1)^{\#S+\#T}\Gamma_{ST}. 
         \label{eq:cist} \end{eqnarray}
\end{theorem}

\textbf{Proof.} By equations (\ref{eq:admd}), (\ref{eq:mbm}), and 
(\ref{eq:bmc}) we have
  $$ A^{-1}= DM^{-1}D^{-1},\ B^{-1}=MDM^{-1}D^{-1},\
     C^{-1}=MDM^{-1}D^{-1}M. $$
Equation (\ref{eq:aist}) is then an immediate consequence of Lemma~
\ref{lemma:minv} and the definition of matrix multiplication.

Since $B^{-1}=MA^{-1}$ we have for fixed $S\supseteq U$ that
  \beas B_{SU}^{-1} & = & \sum_{T\st S\supseteq T\supseteq U} 
    (-1)^{\#T+\#U}A_{TU}\\ 
     & = &  \sum_{T\st S\supseteq T\supseteq U} (-1)^{\#T+\#U}
      \#\{w\in\sn\st \overline{T}\subseteq C(w),\ U\subseteq D(w)\}\\ 
     & = &  \sum_{\overline{T}\st \overline{U}\subseteq
     \overline{T}\subseteq \overline{S}} (-1)^{\#T+\#U} 
      \#\{w\in\sn\st \overline{T}\subseteq C(w),\ U\subseteq D(w)\}.
   \eeas
Equation (\ref{eq:bist}) is now an immediate consequence of the
Principle of Inclusion-Exclusion (or of the equivalence of equations
(\ref{eq:ie1}) and (\ref{eq:ie2})).
Equation (\ref{eq:cist}) is proved analogously to (\ref{eq:bist})
using $C^{-1}=B^{-1}M$. $\ \Box$


\section{Multisets and inversions.}
In this section we consider two further aspects of the connectivity
set: (1) an extension to permutations of a multiset and (2) a
$q$-analogue of Theorem~\ref{thm:main} when the number of inversions
of $w$ is taken into account.

Let $T=\{i_1,\dots,i_k\}_<\subseteq [n-1]$. Define the multiset
  $$ N_T =\{ 1^{i_1},2^{i_2-i_1},\dots,(k+1)^{n-i_k}\}. $$
Let $\fs_{N_T}$ denote the set of all permutations of $N_T$, so
$\#\fs_{N_T} = n!/\eta(T)$; and let $w=a_1 a_2\cdots a_n\in\fs_{N_T}$.
In analogy with equation (\ref{eq:cwdef}) define 
  $$ C(w) = \{i\st a_j<a_k\ \mathrm{for\ all}\ j\leq i<k\}. $$ 
(Note that we could have instead required only $a_j\leq a_k$ rather
than $a_j<a_k$. We will not consider this alternative definition
here.) 

\begin{proposition} \label{prop:ms}
Let $S,T\subseteq [n-1]$. Then
  \beas  \#\{w\in\fs_{N_T}\st C(w)=S\} & = & 
         (\Gamma M)_{\overline{S}\,\overline{T}}\\
              & = & \sum_{U\st U\supseteq \overline{T}}
            \Gamma_{\overline{S}U}\\
          & = & \#\{w\in\sn\st C(w)=S,\ D(w)\supseteq \overline{T}\}.
  \eeas
\end{proposition}

\textbf{Proof.} The equality  of the three expressions on the
right-hand side is clear, so we need only show that
  \beq \#\{w\in\fs_{N_T}\st C(w)=S\}=\#\{w\in\sn\st C(w)=S,\
     D(w)\supseteq \overline{T}\}. \label{eq:multi} \eeq
Let $T=\{i_1,\dots,i_k\}_<\subseteq [n-1]$. Given $w\in\sn$ with
$C(w)=S$ and $D(w)\supseteq\overline{T}$, in $w^{-1}$  replace
$1,2,\dots,i_1$ with 1's, replace $i_1+1, \dots,i_2$ with 2's, etc. It
is easy to check that this yields a bijection between the sets
appearing on the two sides of (\ref{eq:multi}).  $\ \Box$ 

\medskip
Let us now consider $q$-analogues $A(q)$, $B(q)$, $\Gamma(q)$ of the
matrices $A,B,\Gamma$. The $q$-analogue will keep track of the number
inv$(w)$ of inversions of $w=a_1\cdots a_n\in\sn$, where we define
   $$ \mathrm{inv}(w) = \#\{(i,j)\st i<j,\ a_i>a_j\}. $$
Thus define
  $$ \Gamma(q)_{ST} =\sum_{{w\in\sn\atop C(w)=\overline{S},\
        D(w)=T}}q^{\mathrm{inv}(w)}, $$
and similarly for $A(q)_{ST}$ and $B(q)_{ST}$. We will obtain
$q$-analogues of Theorems~\ref{thm:alpha} and \ref{thm:main} with
completely analogous proofs.

Write $\bm{(j)}=1+q+\cdots+q^{j-1}$ and $\bm{(j)!}=\bm{(1)(2)
\cdots (j)}$, the standard $q$-analogues of $j$ and $j!$. Let
$S=\{i_1,\dots,i_k\}_<\subseteq [n-1]$, and define
  $$ \eta(S,q) = \bm{i_1!\,(i_2-i_1)!\cdots (i_k-i_{k-1})!\,
  (n-i_k)!}. $$
Let $T\subseteq [n-1]$, and let $\overline{T}=
\{i_1,\dots,i_k\}_<$. Define 
  $$ z(T) = {i_1\choose 2}+{i_2-i_1\choose 2}+\cdots+
     {n-i_k\choose 2}. $$
Note that $z(T)$ is the least number of inversions of a permutation
$w\in\sn$ with $T\subseteq D(w)$. 

\begin{theorem} \label{thm:qalpha}
We have
  $$ A(q)_{ST}=\left\{ \begin{array}{rl}
        q^{z(T)}\eta(\overline{S},q)/\eta(\overline{T},q), &   
        \mathrm{if}\ \overline{S}\cap T=\emptyset;\\ 0, &
        \mathrm{otherwise}. \end{array} \right. $$
\end{theorem}

\textbf{Proof.} Preserve the notation from the proof of
Theorem~\ref{thm:alpha}. If $(s,t)$ is an inversion of $w$ (i.e.,
$s<t$ and $a_s>a_t$) then for some $0\leq h\leq j$ we have $c_h+1\leq
s<t\leq c_{h+1}$. It is a standard fact of enumerative combinatorics
(e.g., \cite[(21)]{rs:bp}\cite[Prop.~1.3.17]{ec1}) that if
$U=\{u_1,\dots,u_r\}_<\subseteq [m-1]$ then
  \beas \sum_{{v\in\fs_m\atop D(v)\subseteq U}} q^{\mathrm{inv}(v)}
    & = & \displaystyle\bm{\left( {\bm{m}\atop
 \bm{u_1,u_2-u_1,\dots,m-u_r}}\right)}\\ & := &  
 \frac{\bm{(m)!}}{\bm{(u_1)!\,(u_2-u_1)!}\cdots\bm{(m-u_r)!}}, \eeas
a $q$-multinomial coefficient. From this it follows easily that if
 $\overline{U}= \{y_1,\dots,y_s\}_<$ then
  $$ \sum_{{v\in\fs_m\atop D(v)\supseteq U}} q^{\mathrm{inv}(v)} =
    q^{z(T)}\displaystyle\bm{\left( {\bm{m}\atop
     \bm{y_1,y_2-y_1,\dots,m-y_s}}\right)}. $$
Hence we can parallel the proof of
Theorem~\ref{thm:alpha}, except instead of merely counting the number
of choices for the sequence $u=(a_{c_h}, a_{c_h}+1,\dots,a_{c_{h+1}})$
we can weight this choice by $q^{\mathrm{inv}(u)}$. Then
  $$ \sum_u q^{\mathrm{inv}(u)} =    q^{{i_2-i_1\choose
      2}+\cdots+{i_k-i_{k-1}\choose 2}} 
    \displaystyle \bm{\left(
    {\bm{c_{h+1}-c_h}}\atop \bm{i_2-i_1, i_3-i_2,
    \dots,i_k-i_{k-1}}\right)}, $$ 
summed over all choices $u=(a_{c_h}, a_{c_h}+1,\dots,a_{c_{h+1}})$.
Taking the product over all $0\leq h\leq j$ yields 
$q^{z(T)}\eta(\overline{S},q)/\eta(\overline{T},q)$. $\ \Box$

\begin{theorem} \label{thm:qmain}
The matrices $A(q),B(q),\Gamma(q)$ have the following inverses:
  \beas A(q)^{-1}_{ST} & = & (-1)^{\#S+\#T}A(1/q)_{ST}\\
        B(q)^{-1}_{ST} & = & (-1)^{\#S+\#T}\sum_{{ w\in\sn
        \atop  \overline{S}= C(w),\ 
             T\subseteq D(w)}} q^{-\mathrm{inv}(w)}\\
        \Gamma(q)^{-1}_{ST} & = & (-1)^{\#S+\#T}\Gamma(1/q)_{ST}. \eeas
\end{theorem}

\textbf{Proof.} Let $D(q)=(D(q)_{ST})$ be the diagonal matrix with
$D(q)_{SS}= \eta(\overline{S},q)$. Let $Q(q)$ be the diagonal matrix
with $Q(q)_{SS}= q^{z(S)}$. Exactly as for (\ref{eq:admd}),
(\ref{eq:mbm}) and (\ref{eq:bmc}) we obtain 
  \beas A(q) & = & D(q)MD(q)^{-1}Q(q)\\
    M\Gamma(q) M & = & A(q)\\
    B(q) & = & M\Gamma(q) = A(q)M^{-1}. \eeas
The proof now is identical to that of Theorem~\ref{thm:main}. $\ \Box$ 

Let us note that Proposition~\ref{prop:ms} also has a straightforward
$q$-analogue; we omit the details.

\pagebreak

\end{document}